# Que révèle l'activité de validation de démonstration circulaire sur la compréhension de démonstration

Alexis Gautreau[1]

[1]Université Paris Cité, Paris, France; alexisgautreau00@gmail.com

*L'étude présentée s'inscrit dans la recherche sur la validation des démonstration considérée comme un prisme permettant de faire ressortir des phénomènes didactiques liés à la démonstration et à l'activité de démonstration. Nous revisitons le cas déroutant d'élèves de 5e qui valident des démonstrations circulaires. Il s'agit de démonstrations dans lesquelles la conclusion est utilisée comme prémisse d'un pas de déduction. Si ces élèves âgés de 12 à 13 ans identifient correctement l'énoncé final de ces démonstrations comme une reformulation de la conclusion de l'énoncé cible, ils ont toutefois des difficultés à interpréter la manière dont la condition de l'énoncé cible est reprise dans la démonstration. Notre analyse remet en question l'interprétation de ce phénomène par Miyazaki et al., qui attribuent l'acceptation par les élèves de la validité du raisonnement circulaire à une mauvaise compréhension du modus ponens lorsqu'il y en a deux à la suite. Nous proposons une autre explication fondée sur la distinction entre le statut opérationnel et le statut théorique des propositions mathématiques. Cette distinction fournit à la fois une justification théorique de l'invalidité des démonstrations circulaires et un outil didactique pour examiner l'activité argumentative des élèves lorsqu'ils tentent de rejeter de telles démonstrations.*

## Introduction

La recherche présentée est un fragment d'un travail de doctorat dont le but est l'élaboration d'un modèle de compréhension en lecture de démonstrations géométriques (RCGP[1]). Par démonstration, nous entendons raisonnement hypothético-déductif écrit sous forme de texte ; et signalons tout de suite que l'énoncé cible considéré dans ce papier est une implication universellement quantifiée. Notre modèle de compréhension repose à la fois sur des concepts issus de la littérature didactique et sur une analyse des données recueillies lors d'une expérimentation, au cours de laquelle 55 élèves de 12-13 ans ont été confrontés à quatre tâches de validation de démonstration. L'une des démonstrations

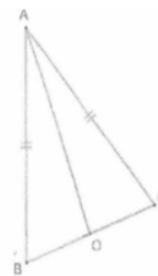

**Claim:** If a triangle is isosceles, then its angles at the base are equal.

**Proof:**
1. We are going to work on triangle $ABC$, which is isosceles at $A$.
2. Let's draw the segment $[AO]$ that divides angle $\widehat{BAC}$ into two equal angles: $\widehat{BAO}$ and $\widehat{CAO}$.
3. We know that $AB = AC$, $\widehat{BAO} = \widehat{CAO}$ and $\widehat{ABO} = \widehat{ACO}$.
4. [Or,] for two triangles to be equal, they only need to have a side of the same length
5. between two equal angles (case of equality ACA).
6. So triangles $ABO$ and $ACO$ are equal.
7. [Or,] by definition, two equal triangles have their homologous angles equal.
8. So angles $\widehat{ABO}$ and $\widehat{ACO}$ are equal.

---

[1] RCGP: abréviation de « Reading Compréhension of Geomtrical Proof » Yang & Lin (Yang & Lin, 2008)

(figure 1) comporte un raisonnement circulaire, c'est-à-dire que la conclusion de la démonstration est mobilisée comme prémisse d'un de ses pas de déduction (ligne 3).

**Figure 1: Démonstration circulaire avec son énoncé cible**

Ce papier se focalise sur ce que l'activité de validation d'un raisonnement circulaire révèle de la compréhension qu'ont les élèves des démonstrations. Nous y présentons le phénomène didactique suivant : certains élèves affirment que la démonstration démontre l'énoncé cible parce que la conclusion de la démonstration correspond à celle de l'énoncé, tout en ne comprenant pas que la démonstration démontre une implication. Plus précisément, ils ne saisissent pas que la démonstration s'ouvre par la mise en hypothèse de la condition de l'énoncé cible.

La première partie expose notre theoretical framework. Celui-ci met l'accent sur l'intérêt didactique de distinguer pour chacune des phrases de la démonstration leur fonction dans le pas de déduction (statut opératoire) de leur position par rapport à la théorie mathématique (statut théorique). La seconde partie s'ouvre sur une analyse critique de la proposition formulée par Miyazaki & al. (2017) concernant l'acceptation d'un raisonnement circulaire. Nous soutenons que cette acceptation résulte d'une incapacité à interpréter le statut d'hypothèse de la phrase de la démonstration qui reprend la condition de l'énoncé cible. Enfin, la troisième partie illustre cette interprétation au moyen d'une étude de cas : l'analyse d'une série d'arguments formulés par un groupe de trois élèves accompagnés d'un enseignant stagiaire placés en position de validateur d'une démonstration circulaire.

## Ancrage théorique

Mejia-Ramos & al. (2012) distinguent dans leur modèle de compréhension de démonstration une composante locale et une composante globale. La composante globale capture la capacité à mettre en œuvre des stratégie de compréhension : résumer la démonstration ; extraire l'idée essentielle… La composante locale correspond à la compréhension des mots, des symboles, des phrases et du rôle des phrases dans le fonctionnement hypothético-déductif. Le contexte d'initiation à la démonstration de notre recherche nous conduit à ne considérer que la composante locale ; nous cherchons à identifier les fonctions relatives à l'organisation hypothético-déductive, que les élèves attribuent aux phrases de la démonstration.

En développant leur *cadre* visant l'évaluation de la compréhension de démonstration, Miyazaki & al. (2017) soulignent que l'outil qu'ils élaborent a pour vocation de modéliser la compréhension de la démonstration en tant qu'« objet structural », qu'ils caractérisent comme « the relational network via deductive reasoning that combines singular and universal propositions. » (Miyazaki et al., 2017, p226). Nous adoptons également cette perspective. Les propositions singulières correspondent sur le plan logique à des propositions non quantifiées, à variables libres : par exemple « We are going to work on the triangle ABC, which is isosceles in A. » ou encore « So angle ABO and angle ACO are equals. » Les propositions universelles qu'évoquent Miyazaki & al. sont les propriétés qui appartiennent au répertoire des propriétés admises ou démontrées sur lequel les élèves s'appuient pour démontrer. Nous nommons *théorie* ce répertoire. Notons que, du point de vue de la logique des prédicats, la théorie contient des énoncés quantifiés universellement mais aussi existentiellement. Cependant, les énoncés de la théorie qui peuplent les démonstrations usuelles dans notre contexte de recherche sont quasi-exclusivement des énoncés universels. Ces derniers s'articulent à des prémisses et une conclusion au sein de ce que nous appelons, à la suite de Duval (1998), *pas ternaire.* Du point de vue de la logique des prédicats, les pas ternaires sont constitués de deux règles d'inférence : l'instanciation universelle (une implication universelle est instanciée en une implication singulière) et le *modus ponens* (combinaison de prémisses singulières et d'une implication singulière, qui se résout en une conclusion singulière).

En plus d'utiliser cette distinction entre proposition singulière et proposition universelle pour modéliser localement la compréhension de démonstration, nous ajoutons que chacune des phrases de la démonstration peut être interprétée doublement : d'une part, selon leur rôle dans le pas ternaire, et, d'autre part, selon leur position par rapport à la théorie. Duval (1998) a stabilisé un lexique utile pour désigner les propositions des démonstrations en fonction de leur rôle dans le pas ternaire, c'est-à-dire de leur *statut opératoire* (Duval, 1998). Une phrase peut être *prémisse* d'un pas en tant que condition d'application d'une propriété universelle ; elle peut également être *énoncé tiers* en tant qu'implication qui noue prémisse et conclusion ; ou encore *conclusion* du pas ternaire. Pour comprendre une démonstration, le lecteur doit interpréter chaque phrase selon son statut opératoire.

Ce n'est toutefois pas suffisant pour modéliser la compréhension de démonstration. Un lecteur peut en effet identifier les statuts opératoires des phrases d'un pas ternaire, sans pour autant comprendre que telle prémisse est mobilisable parce qu'elle est une proposition vraie par supposition, que tel énoncé tiers l'est parce qu'il s'agit d'une implication universellement quantifiée appartenant à la théorie, ou encore que telle conclusion d'un pas ternaire peut être réinvestie en tant que prémisse, du fait qu'elle est vraie par déduction. Nous soutenons dès lors qu'un modèle de compréhension de démonstration doit distinguer nettement, pour chaque phrase, deux statuts : le statut opératoire et le statut théorique.

Le statut théorique d'une phrase indique la position de la phrase par rapport à la théorie. Duval (1998) a uniquement évoqué les statuts théoriques suivants : définitions, axiomes, propriété et conjecture. Le statut théorique d'axiome indique que la vérité de cet énoncé quantifié est assurée par une supposition, tandis que le statut théorique de propriété indique que la vérité de cet énoncé quantifié est garantie par une démonstration. Quant à lui, le statut théorique de conjecture indique que l'énoncé n'est que candidat pour intégrer la théorie. À notre connaissance, Duval n'a pas affecté de statuts théoriques aux propositions singulières présentes dans les démonstrations. Dans le paragraphe suivant, nous présentons de quelle manière nous étendons la signification du concept de statut théorique aux propositions singulières.

Contrairement aux énoncés quantifiés de la théorie qui sont vrais globalement dans la théorie, les propositions singulières d'une démonstration sont vraies dans les limites de la démonstration à laquelle elles appartiennent. Leur vérité est activée par la mise en hypothèse de l'une d'elles en début de démonstration et désactivée en fin de démonstration par l'introduction de l'implication qui constitue la forme instanciée de l'énoncé cible. Les propositions singulières s'ajoutent donc provisoirement à la théorie, le temps de la démonstration. Par ailleurs, leur vérité peut leur être assignée de deux façons : soit les propositions sont vraies par supposition, ce qui signifie que l'auteur ou le lecteur a mis en œuvre en début de démonstration l'opération cognitive suivante : faire l'hypothèse de la condition de l'énoncé cible ; soit ces propositions sont vraies par déduction, ce qui signifie qu'elles sont conclusion d'un pas de déduction au sein de la démonstration. L'assignation de leur vérité intervient donc selon une certaine chronologie hypothético-déductive. En effet, dès la mise en hypothèse en début de démonstration, toutes les propositions singulières de la démonstration sont vraies, même si elles n'ont pas été encore déduites. Cependant, les propositions vraies par déduction ne sont mobilisables que lorsqu'elles sont déjà conclusions d'un pas de déduction. Cette nuance entre *déduite* et *vraie* a des implications didactiques : des élèves considèrent que certaines propositions peuvent être mobilisée en vertu de leur vérité, même si leur déduction intervient ultérieurement. Nous illustrerons plus loin que l'obligation de respecter cette chronologie déductive quand on démontre n'est pas acquise pour tous les élèves. En résumé, nous définissons le statut théorique d'une proposition singulière comme la source d'assignation de la vérité de ces propositions singulières. Ce qui ajoute à la liste des statuts théoriques dressée par Duval deux éléments : *propositions vraies par supposition* et *propositions vraies par déduction*.

Nous nous focalisons dans ce papier sur la compréhension du statut théorique des propositions de la démonstration qui font le lien avec l'énoncé cible. Selden & Selden (1995) dénomment ces phrases *proof framework (épure de preuve)* :

> By a proof framework we mean a representation of the "top-level" logical structure of a proof, which does not depend on detailed knowledge of the relevant mathematical concepts, but which is rich enough to allow the reconstruction of the statement being proved or one equivalent to it. (Selden & Selden, 1995, p. 129)

L'épure de preuve examinée est constitué des deux phrases suivantes :

| 1 | We are going to work on triangle $ABC$, which is isosceles at $A$. |
| 8 | So angles $\widehat{ABO}$ et $\widehat{ACO}$ are equal. |

**Figure 2 : Épure de preuve**

Notons par ailleurs qu'un lecteur peut comprendre le statut théorique de *proposition vraie par supposition* (1$^{re}$ phrase de l'épure de preuve) ou celui de *proposition vraie par déduction* (2$^{e}$ phrase de l'épure de preuve), sans comprendre nécessairement que les phrases de l'épure de preuve constituent respectivement des formes instanciées de la condition et de la conclusion de l'énoncé cible. Dans ce papier, nous ne nous occupons pas de la compréhension de cette instanciation, mais nous prenons en compte, outre les statuts théoriques, le fait que les lecteurs interprètent les deux phrases de l'épure de preuve en tant que reprises-modifications de la condition et de la conclusion de l'énoncé cible. Voici la liste des quatre interprétations attendues de l'épure de preuve qui nous intéressent dans ce texte (table 1).

**Tableau 1 : Interprétations attendues des phrases de l'épure de preuve**

| | |
|---|---|
| Première phrase | Proposition vraie par supposition |
| | Reprise-modification de la condition de l'énoncé cible |
| Seconde phrase | Proposition vraie par déduction |
| | Reprise-modification de la conclusion de l'énoncé cible |

## Analyse de la tâche de validation de démonstration circulaire

Cette partie répond à la question suivante : pourquoi documenter le rapport que des élèves entretiennent aux interprétations attendues de l'épreuve (cf. tableau 1) au moyen de l'activité de validation de démonstration circulaire ?

Nous employons le mot « validation » au sens de Selden & Selden (1995, p. 127) :

> "We use the term validation to describe the process an individual carries out to determine whether a proof is correct and actually proves the particular theorem it claims to prove."

Dans cette définition, valider une démonstration ne se limite pas à vérifier qu'elle est correcte indépendamment de son énoncé cible. Il s'agit également de vérifier qu'elle établit bien le résultat qu'elle prétend démontrer. La démarche de validation, telle que définie par Selden & Selden, porte donc sur le couple énoncé cible / démonstration, et non sur la démonstration détachée de son énoncé cible. Dans notre exemple, où l'énoncé cible est une implication, cela revient à vérifier deux éléments : d'une part, que la condition de l'énoncé cible se retrouve bien dans la démonstration sous la forme d'une hypothèse et que la conclusion de l'énoncé cible correspond effectivement à une proposition déduite en fin de démonstration ; d'autre part, qu'il ne doit pas être fait hypothèse d'une autre proposition que celle correspondant à la condition de l'énoncé cible. Notre démonstration vérifie la première condition mais se heurte à la seconde, puisqu'il est fait hypothèse d'une proposition qui ne correspond pas à la condition de l'énoncé cible.

À notre connaissance, un seul article thématise le raisonnement circulaire dans une perspective didactique. Miyazaki & al. (2017), pour illustrer l'importance d'introduire la règle logique d'instanciation universelle dans un modèle de compréhension de démonstration, analysent un épisode de classe au cours duquel se tient un échange au sujet d'une démonstration circulaire proposée par une élève confrontée à la tâche suivante (figure 3) :

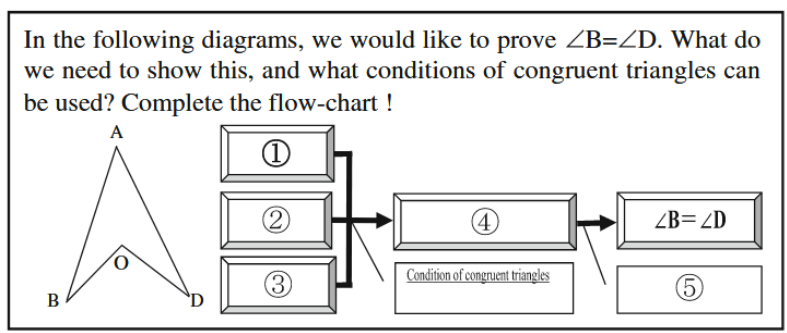

**Figure 3 : Copie de la tâche de Miyazaki & al. (2017)**

Cette démonstration à compléter est présentée sous la forme d'un déductogramme. Celui-ci est conçu de la manière suivante : les rectangles aux contours épais sont destinés à accueillir les propositions singulières, tandis que les rectangles aux contours fins sont réservés aux propriétés quantifiées de la théorie. Notons que cette tâche ne comporte pas d'énoncé cible. Les élèves doivent sélectionner les propositions singulières initiales uniquement à partir de leur connaissance de la figure proposée. La mobilisation de ces propositions singulières qui amorce la démonstration n'est donc pas le fruit de l'opération cognitive de *mise en hypothèse* de la condition de l'implication d'un énoncé cible.

À travers cette tâche, Miyazaki & al. (2017) mettent en évidence le phénomène didactique suivant : l'une des élèves de leur expérimentation maîtrise l'instanciation universelle, ce que les auteurs établissent à l'aide du critère suivant : elle comprend que le pas ternaire repose sur une propriété universellement quantifiée et non sur une de ses formes instanciées. Cependant, cette élève suggère d'inscrire la conclusion de la démonstration (Angle B = Angle D) dans l'un des rectangles placés en entrée du déductogramme – ce rectangle est dédié aux hypothèses, prémisses du premier pas de déduction. Elle maintient ce choix malgré les objections formulées par ses pairs et l'enseignant. Or, cette élève avait réussi antérieurement des tâches analogues, qui ne différaient que par le fait que les démonstrations ne contenaient qu'un seul pas ternaire au lieu de deux. Les auteurs concluent de cette observation que le problème de l'acceptation de la circularité logique provient d'une compréhension insuffisante du *modus ponens*, qui se manifeste lorsque cette règle est appliquée deux fois de suite.

Nous n'interprétons pas cette épisode de classe de la même manière. Pour une raison de logique mathématique d'abord. Miyazaki & al. (2017) considèrent que la démonstration, détachée dans leur exemple de son énoncé cible, est invalide. Plus précisément, ils considèrent que chaque pas de déduction est valide, mais que l'assemblage des deux pas est invalide.

> While this is not problematic when these parts are examined individually, to complete a proof, these two parts […] should be connected by hypothetical syllogism […] Yet in the case of student A this is impossible because the statement to be concluded is used as an assumption, i.e., there is circular reasoning. (Miyazaki et al., 2017, p. 236)

Cette interprétation nous paraît incorrecte. Ce qui noue deux pas de déduction entre eux est le statut théorique de la proposition singulière qui est conclusion du premier *modus ponens*. Ce statut théorique de *proposition vraie par déduction* autorise l'utilisation de cet énoncé en prémisse du second *modus ponens*. En dépit du raisonnement circulaire, la démonstration proposée par l'élève n'est pas invalide. En effet, les propriétés universelles qu'elles mobilisent dans sa démonstration sont bien des propriétés de la théorie, et les propositions singulières sont bien exploitables, puisqu'elles admettent le statut théorique de *proposition vraie par supposition* ou de *proposition vraie par déduction*. Le raisonnement circulaire est donc valide et assigne le vrai à la proposition suivante : « Si les angles B et C sont égaux et *autres conditions*, alors les angles B et C sont égaux ». Ce raisonnement valide, fondé sur une tautologie[2], n'a cependant pas d'intérêt mathématique.

La circularité d'une démonstration ne l'invalide pas, si on la considère indépendamment de son énoncé cible. Invalider un raisonnement circulaire requiert de connaître l'énoncé cible et de montrer que cet énoncé n'est pas ce que démontre la démonstration. Dans notre démonstration, le raisonnement circulaire démontre l'implication suivante : « Si un triangle est isocèle et que ses deux angles à la base égaux, alors ses deux angles à la base sont égaux. » La démonstration de notre tâche, qui n'est pas détachée de son énoncé cible, est donc invalide au sens de Selden & Selden, puisqu'elle ne démontre pas l'énoncé cible. L'épure de preuve présenté en figure 2 correspond bien à l'énoncé cible, mais la proposition « l'angle ABO égal à l'angle ACO » se glisse, en tant qu'hypothèse dans l'épure de preuve, ce qui rend la démonstration invalide.

---

[2] *Tautologie* au sens où l'implication $(A \cap B \rightarrow B)$ est vraie, quelle que soit la valeur de vérité des propositions élémentaires A et B.

L'argument mathématique décisif pour la réfutation du raisonnement circulaire est donc le suivant : la prémisse *Angle ABO égal Angle ACO* du premier pas ternaire, est utilisée comme si elle était supposée vraie, alors que, d'une part, elle n'est pas déduite et, que, d'autre part, elle ne peut être supposée, puisqu'elle n'est pas une forme instanciée de la condition de l'énoncé cible. Les tâches de validation de démonstration circulaire constituent donc de bons candidats pour la mise au jour de la manière dont les élèves interprètent les statuts théoriques des propositions vraies par supposition ainsi que leurs liens avec l'énoncé cible.

## Analyse d'un échange entre élèves autour d'un raisonnement circulaire

Nous analysons ici un épisode de travail en groupe au cours duquel trois élèves (S1, S2 et S3) et un enseignant stagiaire échangent entre eux au sujet d'une démonstration circulaire qu'ils doivent valider. Nous présentons ci-dessous (Tableau 2) le tableau récapitulatif de l'éventuelle adoption ou non des interprétations attendues pour les trois élèves. *P* indique que le lecteur adopte l'interprétation attendue, et *NP* indique que le lecteur n'adopte pas l'interprétation attendue. Les codages indiqués sont justifiés par l'analyse qui suit.

**Table 2: Relation to the four interpretations of S1 and S2**

|  |  | S1 | S2 | S3 |
|---|---|---|---|---|
| 1re phrase | Proposition vraie par supposition | NP | NP | NP |
| | Reprise-modification de la condition de l'énoncé cible | NP | NP | P |
| 2de phrase | Proposition vraie par déduction | NP ? | P | P |
| | Reprise-modification de la conclusion de l'énoncé cible | P | P | P |

Commençons par les interprétations attendues relatives à la seconde phrase de l'épure de preuve. Les trois élèves affirment assez vite que la démonstration établit l'énoncé cible parce que la conclusion de la démonstration affirme l'égalité des angles à la base du triangle isocèle. S2 et S3 comprennent que cette proposition est déduite de la démonstration, c'est moins clair pour S1, qui considère probablement que tout ce qui est vrai (« parce que logique ») puisse être affirmé sans respecter de chronologie déductive.

L'analyse des interprétations de la 1re phrase de l'épure de preuve s'avère plus délicate. Nous la mènerons au moyen d'une synthèse de leurs arguments. S3, perçoit dès la première minute le problème de circularité logique. Elle convainc S2 sans difficulté de l'invalidité de la démonstration en lui signalant que la conclusion de la démonstration est également prémisse d'un pas. Pour S1, la démonstration est bien valide ; cette élève résiste aux différents arguments pendant les 40 minutes de discussion. Les arguments sont de quatre types, que nous exemplifions par des transcriptions d'extraits de l'enregistrement audio de cet épisode. Les deux premiers types d'arguments renvoient à la pertinence du raisonnement en interrogeant la duplication de la conclusion de la démonstration, tandis que les deux derniers renvoient à la validité du raisonnement en interrogeant le droit de mobiliser la prémisse en question.

*(a) L'argument du dédoublement de la conclusion*

S3: Bah oui. Tu veux construire un marteau. Faut heu…T'as un manche, t'as la petite lame et tout, et il faut un clou ! Et pour faire un marteau, tu as besoin d'un autre marteau pour taper. Ça sert à rien de faire un marteau si tu en as déjà un.

S1: Bah si… Tu as créé quelque chose.

L'objection de S1 met en exergue la maladresse de l'analogie : il peut être intéressant de disposer de deux marteaux, même s'il est effectivement inutile de chercher à démontrer la vérité d'un énoncé en le supposant préalablement vrai. Dans tous les groupes étudiés dans notre expérimentation, les enseignants stagiaires, à cours de ressources argumentatives pour convaincre certains élèves de leur groupe, finissent par asséner l'*argument du dédoublement de la conclusion* à la manière d'une jurisprudence à respecter :

Teacher: En mathématiques, tu ne peux pas utiliser quelque chose que tu es en train de démontrer.

De surcroît, cet enseignant, ainsi que les autres enseignants stagiaires participant à l'expérimentation, encouragent en fin de séance les élèves de leur groupe à rédiger, en guise de justification de l'invalidité de la démonstration, des réponses qui relèvent de l'*argument du dédoublement*. Cet argument, qui semble être l'argument privilégié des enseignants et des élèves, ne convainc pas S1 ; elle précise en effet qu'elle ne saisit pas pourquoi il serait superflu de mentionner deux fois la prémisse en question puisqu'elle est vraie.

*(b) L'argument de l'inutilité de ce qui suit la prémisse en question*

Ce second argument la convainc du caractère inutile de la double occurrence de la conclusion.

S3: Alors à ce moment-là, ça veut dire que toi tu dis que la conclu… toute la démonstration elle s'arrête là [Gabrielle pointe la prémisse en question]. T'es d'accord ? Pour toi.

S1: Bah oui.

Cette fois-ci, l'argument n'est plus que la démonstration contient deux occurrences de la conclusion, mais que le texte qui succède à la première occurrence de la proposition est inutile. S1 accepte cet argument, mais en tire la conclusion suivante : le problème de cette démonstration n'est pas la prémisse en question, mais le texte qui suit. Ce dernier est inutile puisque la conclusion de la démonstration est déjà affirmée. Ces objections de S1, que plusieurs élèves de l'expérimentation ont également formulées, pointent la limite des deux premiers types arguments qui s'appuient sur un critère de pertinence, plutôt que sur un critère de validité logique : convaincre de l'inutilité de phrases de la démonstration n'est pas convaincre de l'impossibilité logique d'utiliser une proposition.

*(c) L'argument de l'absence de statut théorique*

S3: [Gabrielle pointe le prémisse erronée] Bah on ne le sait pas… que c'est la conclusion. On ne le sait pas encore.
S1: Bah si on le sait.
S2: Ben non.
S3: Au début, y a rien qui le dit.
S1: Bien sûr que si.
S2: Mais y a rien qui le dit que c'est vrai.
S1: Bien sûr que si.

À la différence de S1, les élèves S3 et S2 ne reconnaissent pas de statut théorique à la prémisse en question, « elle vient de nulle part » ; ils la jugent illégitime. S1 la considère comme évidente, « parce que logique ». L'analyse des autres situations de notre expérimentation suggère que S1 a compris qu'elle ne devait pas se fier à sa perception de la figure ou au mesurage en géométrie déductive, mais que, en revanche, elle cimente en un bloc les propositions *triangle isocèle, côtés égaux et angles égaux,* considérant l'une comme ne pouvant aller sans les autres. Tanguay (2005) a documenté ce phénomène didactique sous le nom de prégnance de la valeur de vérité. Une de nos hypothèses est que cette conception par bloc de l'égalité des côtés et de l'égalité des angles empêche S1 de concevoir que l'une puisse dériver de l'autre selon une certaine chronologie déductive. Ce qui va de pair avec le fait de ne pas discriminer les propositions qui peuvent être mobilisées en tant que prémisse de celles qui ne peuvent l'être. S1 considère donc que cette prémisse a sa place puisqu'il s'agit bien d'un triangle isocèle. Nous ne pouvons pas conclure de cela que S1 ne comprend pas le statut théorique de proposition vraie par supposition ou encore qu'elle ne saisit pas le lien entre la première phrase de l'épure de preuve et la condition de l'énoncé cible. Ses objections à l'argument suivant nous permettrons de conclure.

*(d) L'argument de la comparaison des prémisses selon leur statut théorique*

Pour convaincre S1 de l'absence de raison d'être de la prémisse en question, S3 propose d'analyser le début de la démonstration phrase par phrase afin de donner la raison d'être des deux premiers prémisses et de montrer par contraste l'absence de raison d'être de la troisième prémisse.

| | | |
|---|---|---|
| S1: | | Regarde : « Nous allons travailler dans le triangle isocèle en A. » |
| S3: | | Ça c'est vrai, c'est le triangle ABC, c'est en A qui part des deux côtés. On est d'accord. |
| S1: | | ABC est isocèle en A. Là, déjà, on a le triangle. |
| S3: | | Oui, c'est vrai, la première phrase déjà elle est vraie. Tu peux mettre un petit truc [coche à côté de la phrase], c'est bon. |
| […] | | |
| S1: | | Oui, ben voilà, et après ils mettent là, notre truc [Jeanne pointe la prémisse en question] il est fait. |
| S3: | | Non, et là il nous manque cette info-là. Qu'est-ce qui nous donne cette info-là, que cet angle est égal à cet angle-là ? Personne… |
| S1: | | C'est vrai, regarde : ABO égal ACO, c'est vrai. |

Nous ne discuterons pas ici de la prémisse « Angle BAO = angle CAO » qui demanderait un long développement. Nous ne considérons dans ce papier que la comparaison entre la prémisse *côté égaux* et la prémisse erronée. La prémisse *Côté égaux* est acceptée par S3, S2 et S1 sans débats, sous l'argument que le triangle est isocèle. En revanche, le statut théorique de la première phrase *Triangle ABC isocèle en A* n'est jamais explicité comme une hypothèse ou une reprise-modification de l'énoncé cible par S3 et S2. S1 insistera d'ailleurs sur le fait que si l'on peut affirmer cette proposition, il n'y a pas de raison de ne pas affirmer aussi la prémisse erronée, ce qui prouve qu'elle ne fait pas le lien entre condition de l'énoncé cible et première phrase de l'épure de preuve.

Des données extraites d'autres situations prouvent que S3 comprend la nécessité d'une reprise-modification de la condition de l'énoncé cible. En revanche, il semble que S2, quoiqu'il juge illégitime la prémisse en question depuis le début de la discussion, parce que « y a rien qui le dit que c'est vrai », ne tienne pas compte du lien entre la première phrase de l'épure de preuve et la condition de l'énoncé cible. Il n'interprète pas non plus le statut théorique de proposition vraie par supposition, comme nous le suggère l'indice suivant : à la fin du travail en groupe, S2 propose, dans le but de corriger la démonstration, de remplacer « $\widehat{ABO} = \widehat{ACO}$ » par « *Je pense* que ABO égal ACO ». Cette proposition suggère que, pour S2, le problème n'est pas la position de la proposition dans la démonstration, mais la certitude avec laquelle elle est affirmée. Il tente donc de la rendre mobilisable en modalisant sa vérité.

**Conclusion**

Ce papier pointe une difficulté répandue, rencontrée par un groupe d'élèves de notre expérimentation. Tout en étant en mesure de faire le lien entre la conclusion de l'énoncé cible et la dernière phrase de l'épure de preuve, certains élèves ne sont pas en mesure, d'une part, de distinguer les prémisses mobilisables des autres, et, d'autre part, d'exprimer, même maladroitement, le principe de mise en hypothèse de la condition de l'énoncé cible. Les discussions entre élèves au sujet de l'activité de validation de démonstration circulaire ont constitué un levier important pour mettre au jour ce phénomène didactique. Ce dernier suggère de se pencher sur la question de la mise en œuvre d'une étude exploratoire sur les formes que pourrait prendre un enseignement du statut théorique de *proposition déduite par supposition* et du lien que des propositions de ce type nouent avec l'énoncé cible.

**Références**